\documentclass[11pt]{article}

\usepackage{amsmath,amssymb,url}

\newtheorem{theorem}{Theorem}
\newtheorem{ass}{Theorem}[section]
\newtheorem{prop}[ass]{Proposition}
\newtheorem{lemma}[ass]{Lemma}
\newtheorem{conj}{Conjecture}
\newtheorem{definition}[ass]{Definition}

\newcommand{\qed}{\hspace*{\fill} \rule{7pt}{7pt}\vspace{5pt}}
\newcommand{\Proof}{\noindent{\bf Proof.}\ \ }

\newcommand{\F}{\mathcal F}
\newcommand{\E}{\mathcal E}
\newcommand{\I}{\mathcal I}
\renewcommand{\H}{\mathcal H}

\begin{document}

\title{On the exact maximum induced density of almost all graphs and their inducibility}

\author{
Raphael Yuster
\thanks{Department of Mathematics, University of Haifa, Haifa
31905, Israel. Email: raphy@math.haifa.ac.il. This research was supported in part by the Israel Science Foundation
(grant No. 1082/16).}
}

\date{}

\maketitle

\setcounter{page}{1}

\begin{abstract}

Let $H$ be a graph on $h$ vertices. The number of induced copies of $H$ in a graph $G$ is denoted by $i_H(G)$.
Let $i_H(n)$ denote the maximum of $i_H(G)$ taken over all graphs $G$ with $n$ vertices.

Let $f(n,h) = \Pi_{i}^h a_i$ where $\sum_{i=1}^h a_i = n$ and the $a_i$ are as equal as possible. Let $g(n,h) = f(n,h) + \sum_{i=1}^h g(a_i,h)$.
It is proved that for almost all graphs $H$ on $h$ vertices it holds that $i_H(n)=g(n,h)$ for all $n \le 2^{\sqrt{h}}$.
More precisely, we define an explicit graph property ${\cal P}_h$ which, when satisfied by $H$, guarantees that $i_H(n)=g(n,h)$ for all $n \le 2^{\sqrt{h}}$.
It is proved, in particular, that a random graph on $h$ vertices satisfies ${\cal P}_h$ with probability $1-o_h(1)$.
Furthermore, all extremal $n$-vertex graphs yielding $i_H(n)$ in the aforementioned range are determined.

We also prove a stability result. For $H \in {\cal P}_h$ and a graph $G$ with $n \le 2^{\sqrt{h}}$ vertices satisfying $i_H(G) \ge f(n,h)$,
it must be that $G$ is obtained from a balanced blowup of $H$ by adding some edges inside the blowup parts.

The {\em inducibility} of $H$ is $i_H = \lim_{n \rightarrow \infty} i_H(n)/\binom{n}{h}$.
It is known that $i_H \ge h!/(h^h-h)$ for all graphs $H$ and that a random graph $H$ satisfies almost surely that
$i_H \le h^{3\log h}h!/(h^h-h)$. We improve upon this upper bound almost matching the lower bound.
It is shown that a graph $H$ which satisfies ${\cal P}_h$ has $i_H =(1+O(h^{-h^{1/3}}))h!/(h^h-h)$.

\end{abstract}

\section{Introduction}

We consider the basic problem of maximizing the number of induced copies of a given graph $H$ in a graph with $n$ vertices.
This problem has been studied already for four decades by various researchers such as in
\cite{BHLP-2016,BEHJ-1995,BNT-1986,BS-1994,EL-2015,exoo-1986,FV-2012,HHN-2014,ht-2017,hirst-2014,huang-2014,PG-1975}
and its answer (even asymptotically) seems difficult already for some very small graphs.
Our main result gives an exact answer to this problem for almost all graphs $H$, and for all $n$ that are not too huge.

Let $i_H(G)$ denote the number of induced copies of $H$ in a graph $G$.
The related extremal graph-theoretic parameter of interest is therefore the maximum of $i_H(G)$ taken over all graphs $G$ with $n$ vertices,
denoted by $i_H(n)$. It is a simple exercise to prove that the maximum induced density of $H$ in graphs with $n$ vertices,
namely $i_H(n)/\binom{n}{h}$ tends to a limit,
which is denoted by $i_H = \lim_{n \rightarrow \infty}i_H(n)/\binom{n}{h}$ and is called the {\em inducibility} of $H$.

The first to define and study these natural parameters were Pippenger and Golumbic \cite{PG-1975}.
Let $f(n,h) = \Pi_{i}^h a_i$ where $\sum_{i=1}^h a_i = n$ and the $a_i$ are as equal as possible.
Let $g(n,h)$ be the function recursively defined as $g(n,h)=0$ if $n < h$ and otherwise by $g(n,h) = f(n,h) + \sum_{i=1}^h g(a_i,h)$.
Having noticed that the so-called {\em nested blowups} of $H$ on $n$ vertices (see definition in the next section)
contain at least $g(n,h)$ induced copies of $H$, namely that $i_H(n) \ge g(n,h)$, they used this fact to prove that
$i_H \ge h!/(h^h-h)$ for every graph $H$. However, as it turns out, determining $i_H$ and, moreover, determining $i_H(n)$ seems to be a very difficult task in general.
In fact, there are only very few families of graphs and sporadic cases for which these are known.

It is easy to see that $i_H(n)=i_{H^C}(n)$ where $H^C$ is the complement of $H$. Clearly, $i_{K_h}(n)=i_{I_h}(n)=\binom{n}{h}$
and that cliques and their complements are the only graphs with unit inducibility.
Applying a result of Goodman \cite{goodman-1959} it is also a simple exercise to determine $i_{K_{1,2}}(n)$ and deduce that $i_{K_{1,2}}=3/4$
thereby obtaining a complete solution for all graphs on three vertices.
The inducibility of most graphs on four vertices is known, see \cite{EL-2015,exoo-1986}. However, the inducibility of $P_4$, the path on four vertices,
is still not completely determined. See the website of the flag algebra software of Vaughan \cite{flagmatic-site} for the best upper bound obtained using the flag algebra method of Razborov \cite{razborov-2007}
and Even-Zohar and Linial \cite{EL-2015} for the best lower bound. More generally, the inducibility and maximum induced density of complete bipartite graphs is well understood
\cite{BNT-1986,BS-1994,PG-1975}. For example, $i_{K_{h,h}}=\binom{2h}{h}/4^h$.
Bollob{\'a}s, Egawa, Harris, and Jin \cite{BEHJ-1995} proved that if $H$ is a sufficiently large balanced blowup of some complete graph $K_r$,
then $i_H(n)$ is obtained by blowups of $K_r$. Recently, Hatami, Hirst, and Norine \cite{HHN-2014} proved that
if $H$ is a sufficiently large balanced blowup of some graph $K$, then any graph which attains $i_H(n)$ for $n$ sufficiently large
must itself essentially be a blow-up of $K$.

An intriguing conjecture raised in the paper of Pippenger and Golumbic \cite{PG-1975}, which is yet unsolved, is that the inducibility of
the cycle $C_h$ for $h \ge 5$ is equal to the aforementioned lower bound $h!/(h^h-h)$. The case $h=5$ was only very recently solved by
Balogh, Hu, Lidick{\`y}, and Pfender \cite{BHLP-2016} with sophisticated application of flag algebra. In fact, they have proved that
$i_{C_5}(5^k)=g(5^k,5)$ and that the unique extremal graph is the corresponding nested blowup of $C_5$.
They also proved that for $n$ sufficiently large, $i_{C_5}(n) = g(n,5)$.
So, in this case nested blowups (as opposed to balanced blowups) are extremal graphs.
Nested blowups appear as extremal constructions also in many other extremal problems, as shown by Pikhurko \cite{pikhurko-2014}.

Which other graphs have the property that their inducibility is close (or equal) to the generic lower bound $h!/(h^h-h)$?
Clearly, if we can prove that the nested blowups are extremal for $c_H(n)$ for all $n$, then a by-product would be that $i_H=h!/(h^h-h)$.
Likewise, if we can prove that the nested blowups are extremal for $c_H(n)$ for a very large $n$ as a function of $H$, then we can get
very close to $h!/(h^h-h)$. The main result of this paper does the latter, and does it for almost all graphs $H$.

Recall a random graph on $h$ vertices is the distribution $G(h,1/2)$ on $h$-vertex graphs where each pair of vertices is an edge with probability $1/2$.
Recall also that a graph property ${\cal P}_h$ is a family of $h$-vertex graphs closed under isomorphism.
We say that a graph $H \sim G(h,1/2)$ satisfies ${\cal P}_h$ asymptotically almost surely
(alternatively, that ${\cal P}_h$ is satisfied by almost all graphs), if $\Pr [H \in {\cal P}_h] = 1-o_h(1)$.
In Section 3, we describe an explicit graph property ${\cal P}_h$ that we call {\em strongly asymmetric graphs}.
The following two theorems together form the crux of our result.

\begin{theorem}\label{t:1}
Let $H \sim G(h,1/2)$.
$$
\Pr [H \in {\cal P}_h] = 1-o_h(1)\;.
$$
\end{theorem}

\begin{theorem}\label{t:main}
$ $\begin{enumerate}
\item
For all $n \le 2^{\sqrt{h}}$, if $H \in {\cal P}_h$, then $i_H(n)=g(n,h)$.
Furthermore, the extremal graphs yielding $i_H(n)$ are precisely all the nested blowups of $H$
on $n$ vertices.
\item
For all $n \le 2^{\sqrt{h}}$, if $H \in {\cal P}_h$ and $G$ is an $n$-vertex graph with $i_H(G) \ge f(n,h)$, then
$G$ must be a balanced blowup of $H$ with some edges added inside the blowup parts.
\end{enumerate}
\end{theorem}
Notice that Theorem \ref{t:main} is a completely deterministic statement. We emphasize this since there are graphs that satisfy ${\cal P}_h$ and that are
quite far from a typical element of $G(h,1/2)$ (for example, there are graphs in ${\cal P}_h$ that contain linear sized cliques).
Another point to mention is that one can extend the range for which Theorem \ref{t:main} holds to all $n \le 2^{h^{1-\epsilon}}$
for a given $\epsilon > 0$, at the price of redefining ${\cal P}_h = \emptyset$ when $h$ is sufficiently small, depending on $\epsilon$.

It is known that random graphs $H$ on $h$ vertices satisfy $i_H \le h^{O(\log h)}h!/(h^h-h)$ asymptotically almost surely, as has been observed by Even-Zohar and Linial. Using our result we can improve this upper bound by essentially removing the $h^{O(\log h)}$ factor, thereby coming extremely close to the
generic lower bound.
\begin{theorem}\label{t:inducibility}
If $H \in {\cal P}_h$, then
$$
i_H \le \frac{h!}{h^h-h} \cdot \left( 1+\frac{4}{h^{h^{1/3}}} \right)\;.
$$
\end{theorem}

The rest of this paper has the following structure. The next section consists of definitions and description of the
basic objects required for the rest of the paper. Strongly asymmetric graphs are defined in Section 3, where we also prove that
almost all graphs are strongly asymmetric, namely Theorem \ref{t:1}. The proof of the main result, Theorem
\ref{t:main}, appears in Section 4. Its main ingredient uses a concept of {\em consistent large sets} defined there.
The proofs of the properties of consistent large sets and the main lemma showing that they can be grouped to a single consistent large set are given in
Sections 5 and 6 respectively. Section 7 shows how to apply our main result to inducibility, proving Theorem \ref{t:inducibility}.
The final section contains a conjecture and some concluding remarks.

\section{Preliminaries}

Throughout this paper we generally do not omit floors and ceilings as these are important for the proof to work for the case
where $n$ is close to $h$, as we do not want to treat this (seemingly easier) case separately.
For positive integers $h$ and $n$, a sequence of nonnegative integers $a_1,\ldots,a_h$ whose sum is $n$ is an {\em $(h,n)$-partition}.
We say that the sequence is {\em equitable} if any two elements in the sequence differ by at most $1$.
Otherwise, we say that the sequence is non-equitable.
Clearly, if $n \pmod h = t$ where $0 \le t \le h-1$, then an equitable $(h,n)$-partition
has $t$ elements of order $\lceil n/h \rceil$ and $h-t$ elements of order $\lfloor n/h \rfloor$.

Let $f(n,h)$ be the product of the elements of an equitable $(h,n)$-partition. Then,
$$
f(n,h) = \lceil n/h \rceil^t \lfloor n/h \rfloor^{h-t}\;.
$$
We note that $f(n,h)=0$ for $n < h$ and that $f(h,h)=1$. More generally, $f(kh,h)=k^h$ for a positive integer $k$.
The following lemma summarizes a few properties of $f(n,h)$ and non-equitable partitions.
\begin{lemma}\label{l:simple}
$$
f(n-1,h) = \frac{f(n,h)}{\lceil n/h \rceil}\left( \lceil n/h \rceil -1 \right)\;.
$$
Consequently,
$$
f(n-k,h) \le f(n,h)\left(1-\frac{1}{\lceil n/h \rceil} \right)^k\;.
$$
If $a_1,\ldots,a_h$ is a non-equitable $(h,n)$-partition, then
$$
\Pi_{i=1}^h a_i \le f(n,h)\left(1-\frac{1}{\lceil n/h \rceil^2}\right)\;.
$$
\qed
\end{lemma}
\Proof
The first part, as well as its obvious consequence are straightforward hence we prove only the third part.
Consider placing $n$ elements in $h$ bins where the $i$th bin contains $a_i$ elements.
One can obtain such a placement by starting with an equitable $(h,n)$-partition of the elements into the bins
and repeatedly moving elements from some bin whose current size is $a$ to a bin whose current size is $b$
where $a \le b$. After each move the product of the bin sizes reduces by a factor of $(a-1)(b+1)/ab \le 1-1/(ab)$.
Before the first move we have $a \le b \le \lceil n/h \rceil$. So already after the first move the product reduces by a factor
of at most $1-\frac{1}{\lceil n/h \rceil^2}$. \qed

We next define a recursive variant of $f(n,h)$.
Let $h$ and $n$ be positive integers. Let $g(n,h)=0$ for $n < h$. For $n \ge h$, define
$g(n,h) = f(n,h) + \sum_{i=1}^h g(a_i,h)$ where $a_1,\ldots,a_h$ is an equitable $(h,n)$-partition.
It is immediately observed that $g(n,h)=f(n,h)$ for all $n \le h(h-1)$ while $g(n,h) > f(n,h)$ for $n > h(h-1)$.
When $n=h^k$ for a positive integer $k$ we obtain by induction that
\begin{equation}\label{e:ghk}
g(h^k,h)= \sum_{i=1}^k h^{(k-1)h-(i-1)(h-1)} = \frac{h^{h(k-1)}(1- h^{k(1-h)})}{1-h^{1-h}}\;.
\end{equation}

Let $H$ be a graph on $h$ vertices and assume that $V(H)=[h]=\{1,\ldots,h\}$.
Let $a_1,\ldots,a_h$ be an $(h,n)$-partition, where the $a_i$ are positive integers. The graph $H(a_1,\ldots,a_h)$ is defined as follows.
Its vertex set is the disjoint union of independent sets $A_1,\ldots,A_h$ with $|A_i|=a_i$, hence it has $n$ vertices.
Its edges are defined as follows. For each pair $i,j$ of distinct vertices of $H$, if $ij \in E(H)$
then the bipartite graph in $H(a_1,\ldots,a_h)$ between $A_i$ and $A_j$ is complete, whereas if $ij \notin E(H)$
then the bipartite graph in $H(a_1,\ldots,a_h)$ between $A_i$ and $A_j$ is empty. We say that $H(a_1,\ldots,a_h)$ is a {\em blowup} of $H$.
We call $A_i$ the {\em part} of the blowup corresponding to vertex $i \in V(H)$.
We call a pair of distinct vertices of the blowup a {\em part pair} if both are in the same part
and call them a {\em blowup pair} if they are in distinct parts.
If $a_1,\ldots,a_h$ is an equitable $(h,n)$-partition, then $H(a_1,\ldots,a_h)$ is called an $n$-vertex {\em balanced blowup} of $H$
and is denoted by $H(n)$.

A few simple facts regrading blowups follow. An induced copy of $H$ in $H(a_1,\ldots,a_h)$ can be obtained by selecting one vertex from each part.
This immediately shows that $i_{H}(H(a_1,\ldots,a_h)) \ge \Pi_{i=1}^h a_i$.
In particular, an $n$-vertex balanced blowup of $H$ gives that $i_{H}(H(n)) \ge f(n,h)$ which implies that $i_{H}(n) \ge f(n,h)$.
Notice, however that the last inequality is not sharp for all $n > h(h-1)$. Indeed, consider  $H(n)$.
One of the parts has order at least $h$. Now, instead of letting this part be an independent set, add an induced copy of $H$ to it.
This adds at least one additional induced copy of $H$ to the graph, while keeping all the induced copies having one vertex in each part.
Hence, $i_{H}(n) > f(n,h)$ for $n > h(h-1)$.

The last example immediately triggers the following construction, so called the family of $n$-vertex {\em nested balanced blowups} of $H$,
denoted by $H^*(n)$. For $n < h$ we define $H^*(n)$ to be the set of all graphs on $n$ vertices.
For $n \ge h$ we define $H^*(n)$ as follows.
Take a balanced blowup $H(n)$. Now, replace each independent set $A_i$ with an element of $H^*(a_i)$ (note: if $a_i=a_j$ we
are allowed to replace $A_i$ with an element of $H^*(a_i)$ that is different from the element of $H^*(a_i)$ that replaced $a_j$).
An element of $H^*(n)$ is called a nested balanced blowup of $H$.

A few simple facts regrading nested balanced blowups follow from the definition.
First notice that for all $h \le n \le h(h-1)$ we have that the elements of $H^*(n)$ are just the usual balanced blowups of $H$.
Another obvious but interesting point to observe is that for a positive integer $k$, if $n=h^k$, then $H^*(n)$ has only one element.
Also, by immediate induction, $i_H(X) \ge g(n,h)$ for any $X \in H^*(n)$ which implies that $i_{H}(n) \ge g(n,h)$.
Notice, however, that this inequality is not always tight. In fact, in some cases, $H^*(n)$ can itself contain an element having more than $g(n,h)$ induced copies of $H$.
Consider for example, $H=K_{1,2}$. It is immediate to check that the unique element of $K_{1,2}^*(9)$ contains more than $3^3+3=g(9,3)$ induced copies of $K_{1,2}$.

We end this section by recalling the definition of isomorphism between graphs.
Let $H_1$ and $H_2$ be two graphs of the same order. A bijective function $f:V(H_1) \rightarrow V(H_2)$ is called an {\em isomorphism between $H_1$ and $H_2$} if
$(u,v) \in E(H_1)$ if and only if $(f(u),f(v)) \in E(H_2)$. We say that $H_1$ and $H_2$ are {\em isomorphic} if there exists an isomorphism between them.
An isomorphism between a graph and itself is called an {\em automorphism}. If $f$ is a bijection between $H_1$ and $H_2$,
then a stationary point of $f$ is a vertex $v$ such that $f(v)=v$. Notice that the number of stationary points is at most $|V(H_1) \cap V(H_2)|$.

\section{Strongly asymmetric graphs}

The purpose of this section is to define the family of strongly asymmetric graphs and prove that almost all graphs
reside in this family. However, before we do that, we need a couple of definitions.
In what follows, {\em changing} an edge in a graph means either removing it, or adding it. 

\begin{definition}[{\bf blowup $m$-far; blowup $m$-close}]
Let $H_1$ and $H_2$ be two graphs of the same order, where $H_1$ is a blowup of some other graph $H$.
We say that $H_1$ and $H_2$ are {\em blowup $m$-close} if one can change at most $m$ blowup edges of $H_1$ and arbitrarily change
part edges of $H_1$ such that after the change, the obtained graph is isomorphic to $H_2$. We say that they are {\em blowup $m$-far} if they are not blowup $m$-close.
\end{definition}
Notice that if $H_1$ is a trivial blowup (i.e. $H_1$ is just an arbitrary graph), then the latter definition just means that
we can change at most $m$ edges of $H_1$ such that after the change, the obtained graph is isomorphic to $H_2$.
In this case we say that the graphs are $m$-close or $m$-far. More specifically, we need the following definition:
\begin{definition}[{\bf $m$-far isomorphism; $m$-close isomorphism}]
Let $H_1$ and $H_2$ be two graphs of the same order, and let $f$ be a bijection from $V(H_1)$ to $V(H_2)$.
We say that $f$ is {\em $m$-close to an isomorphism between $H_1$ and $H_2$} if one can change at most $m$ edges of $H_1$ such that $f$
is an isomorphism between the obtained graph and $H_2$. Otherwise, we say that $f$ is {\em $m$-far from an isomorphism between $H_1$ and $H_2$}.
\end{definition}

\begin{definition}[{\bf agreement set; distinguishing set}]
Let $H$ be a graph and let $S \subset V(H)$. Let $A(S) \subseteq V(H) \setminus S$ be the set of vertices with the following property.
Each $w \in A(S)$ is either adjacent to all vertices of $S$ or non-adjacent to all vertices of $S$.
We call $A(S)$ the {\em agreement set} of $S$.
We say that $Q \subset V(H)$ is a {\em distinguishing set} if $Q$ is not contained in the agreement set of any two vertices outside of $Q$.
\end{definition}

Let $H$ be a graph on vertex set $[h]$.
Recall from the previous section that a permutation $\pi : [h] \rightarrow [h]$ is an automorphism of $H$
if $(\pi(i), \pi(j))$ is an edge of $H$ if and only if $(i,j)$ is an edge of $H$.
Clearly, the set of all automorphisms of $H$ is a group with respect to the composition operator, denoted by $Aut(H)$.
Following Erd\"os and R\'enyi \cite{ER-1963}, we say that $H$ is {\em asymmetric} if $Aut(H)$ consists only of the identity permutation
and {\em symmetric} otherwise.

Is is easy to verify that all graphs on at most $5$ vertices (and more than one vertex) are symmetric.
However, there are already asymmetric graphs on $6$ vertices.
The smallest one (with respect to the number of edges) 
is obtained from the path on vertices $1,2,3,4,5$ (in this order) by adding vertex $6$
and connecting it to vertices $3$ and $4$.
Erd\"os and R\'enyi \cite{ER-1963} proved that almost all graphs are asymmetric.

We will need a stronger notion of asymmetry in graphs.
For a positive integer $h$ we say that a graph $H$ on $h$ vertices is {\em strongly asymmetric},
or, synonymously, that it belongs to the family ${\cal P}_{h}$, if all the following conditions hold for $H$:
\begin{enumerate}
\item
The degree of every vertex of $H$ is larger than $0.4h$ and smaller than $0.6h$.
\item
For every pair of vertices, the order of their agreement set is not larger than $0.55h$
and for any triple of vertices, the order of their agreement set is not larger than $0.3h$.
\item
$H$ has a distinguishing set of order at most $3\log h$ \footnote{Throughout this paper all logarithms are in base $2$.}.
\item
Let $S \subset [h]$ with $|S| \le 0.7h$.
Let $B$ be any blowup of $H[S]$ with at least $0.8h$ vertices and at most $h$ vertices, and let $K \subseteq [h]$.
Furthermore, assume that each blowup part of $B$ contains no more than $h/100$ vertices.
Then $B$ is blowup $10^{-5}h^2$-far from $H[K]$.
\item
Let $J \subset [h]$ and $K \subset [h]$ where $|J|=|K|=\lceil 0.7h \rceil$
and let $\pi$ be a bijection from $J$ to $K$ having at least $0.1h$ non-stationary points.
Then $\pi$ is $10^{-5}h^2$-far from an isomorphism between $H[J]$ and $H[K]$.
\end{enumerate}
It should be noted that although not stated explicitly, strong asymmetry implies asymmetry. Indeed, suppose $H$ is strongly asymmetric.
Then, Condition 5 implies that any automorphism of $H$ has at most $0.1h$-non stationary points.
But if an automorphism is non-trivial this means that there is a pair of distinct vertices whose agreement set
is at least $0.9h-2$, which is impossible by Condition 2.

We next show that almost all graphs on $h$ vertices do in fact belong to ${\cal P}_h$.
Recall that ${\cal G}(h,\frac{1}{2})$ is the probability space of all graphs on vertex set $[h]$
where each pair of vertices are connected with an edge with probability
$\frac{1}{2}$, and the $\binom{h}{2}$ choices are independent.
\begin{theorem}\label{t:phe}
For all $h$ sufficiently large, if $H \sim G(h,1/2)$, then
$$
\Pr [H \in {\cal P}_h] \ge 1-\frac{2}{h}\;.
$$
\end{theorem}
Hence Theorem \ref{t:1} is just a more explicit restatement of Theorem \ref{t:phe}.
The proof of Theorem \ref{t:phe} follows from the five following lemmas, each considering one item in the definition of ${\cal P}_{h}$.
While the first three are very elementary and their proofs are only presented for completeness, the latter two are more technical.
\begin{lemma}\label{l:cond-1}
If $h$ is sufficiently large and $H \sim G(h,1/2)$, then with probability at least $1-1/(5h)$ all degrees of $H$ are larger than $0.4h$ and smaller than $0.6h$.
\end{lemma}
\Proof
The degree of a vertex of $H$ is a random variable with distribution $Bin(h-1,1/2)$ so its expectation is $(h-1)/2$. The probability that it deviates
from its expectation by more than a constant factor is exponentially small in $h$, and there are only $h$ vertices to consider.
\qed
\begin{lemma}\label{l:cond-2}
If $h$ is sufficiently large and $H \sim G(h,1/2)$, then with probability at least $1-1/(5h)$, $H$ has the property that for every pair of vertices, the order of their agreement set is not larger than $0.55h$
and for any triple of vertices, the order of their agreement set is not larger than $0.3h$.
\end{lemma}
\Proof
Consider a pair of vertices. Their agreement set is a random variable with distribution $Bin(h-2,1/2)$
so its expectation is $(h-2)/2$. The probability that it deviates from its expectation by more than a constant factor is exponentially small in $h$, and there are less than $h^2$
pairs to consider. Consider a triple of vertices. Then their agreement set is a random variable with distribution $Bin(h-3,1/4)$
so its expectation is $(h-3)/4$. The probability that it deviates from its expectation by more than a constant factor is exponentially small in $h$, and there are less than $h^3$
triples to consider.
\qed
\begin{lemma}\label{l:cond-3}
If $h$ is sufficiently large and $H \sim G(h,1/2)$, then with probability at least $1-1/h$,
$H$ has a distinguishing set of at most $\lfloor 3 \log h \rfloor$ vertices.
\end{lemma}
\Proof
Fix some $Q \subset [h]$ with $|Q|=\lfloor 3 \log h \rfloor$.
We prove that with probability at least $1-1/h$, $Q$ is a distinguishing set.
For two given vertices outside of $Q$, the probability that $Q$ is contained in their agreement set is precisely $2^{|-Q|}$.
As there are less than $h^2/2$ pairs of vertices outside of $Q$, the probability that $|Q|$ is not distinguishing is (by the union bound)
at most $(h^2/2) 2^{|-Q|}$ which satisfies
$$
\frac{h^2}{2}2^{|-Q|} \le \frac{h^2}{2}2^{1-3 \log h} = \frac{1}{h}\;.
$$
\qed
\begin{lemma}\label{l:cond-4}
If $h$ is sufficiently large and $H \sim G(h,1/2)$, then the following holds.
Let $S \subset [h]$ with $|S|=s \le 0.7h$. Let $(a_1,\ldots,a_s)$ be an $(s,k)$-partition with $h \ge k \ge 0.8h$ and with $a_i \le h/100$ for
$i=1,\ldots,s$. Let $K \subseteq [h]$ with $|K|=k$. Then, the probability that $B=H[S](a_1,\ldots,a_s)$ is blowup $10^{-5}h^2$-close to $H[K]$
is at most $e^{-0.001h^2}$.
\end{lemma}
\Proof
Denote the parts of $B$ by $A_1,\ldots,A_s$ where $|A_i|=a_i \le h/100$.
Let us fix a set $M$ of $m \le \lfloor 10^{-5}h^2 \rfloor$ blowup pairs of $B$ and fix a bijection $\pi\,:\, V[B] \rightarrow K$.
Let $B_M$ be the graph obtained from $B$ after changing $M$ (making pairs in $M$ that are edges into non-edges and making
pairs in $M$ that are non-edges into edges).
Call a vertex of $B$ {\em problematic} if it appears in at least $h/30$ elements of $M$.
Notice that being problematic has nothing to do with the structure of $H$, since it only depends on $S$ and $M$ ($S$ is just a subset of $[h]$
and $M$ is just a subset of blowup pairs, while blowup pairs only depend on the partition).
As each pair in $M$ contributes
two to the count towards being problematic, the number of problematic vertices is at most
\begin{equation}\label{e:num-prob}
\frac{2m}{h/30} \le 0.0006h\;.
\end{equation}

Consider now the following maximum set of pairs $\{x_1,y_1\},\ldots,\{x_\ell,y_\ell\}$ that are pairwise disjoint (i.e altogether they contain $2\ell$ vertices)
and have the following property.
For all $j=1,\ldots,\ell$, $\{x_j,y_j\}$ both belong to the same part (no matter which) and both are non problematic.
We lower bound $\ell$ as follows.
Since $\ell$ is maximum, from each part $A_i$ for $i=1,\ldots,s$, there is at most one uncovered non-problematic vertex in our set of pairs.
Thus,
\begin{equation}\label{e:l}
\ell \ge \frac{k-0.0006h-s}{2} \ge 0.04h
\end{equation}
where we have used that $k \ge 0.8h$ and $s \le 0.7h$.
Consider some $\{x_j,y_j\}$ pair and suppose they belong to the part $A_i$.
Although both $x_j$ and $y_j$ are non-problematic, they can still occur each in at most $h/30$ elements of $M$,
so after ignoring the at most $2\cdot h/30=h/15$ vertices in pairs in $M$ to which at least one of $x_j$ or $y_j$ belongs
and after ignoring all the problematic vertices
we still have a set $Z$ of vertices of $B$ of size at least
$$
k-0.0006h-\frac{h}{15}-a_i \ge 0.8h-0.006h-\frac{h}{15}-\frac{h}{100} \ge 0.7h
$$
that has the following property: For each $z \in Z$, $(x_j,z)$ and $(y_j,z)$ are blowup pairs that are not in $M$.
Again, we stress that $Z$ has nothing to do with the structure of $H$ since it only depends on $S$, on $M$, and on problematic vertices.
In other words, regardless of the structure of $H$, in the graph $B_M$, the set $Z$ will be in the agreement set of $x_j$ and $y_j$.
Now consider $H[K]$. In order for $\pi$ to be an isomorphism between $B_M$ (regardless of any additional changes to part pairs of $B$)
and $H[K]$ we need for the agreement set of each pair
$\{\pi(x_1),\pi(y_1)\},\ldots,\{\pi(x_\ell),\pi(y_\ell)\}$ to be of size at least $0.7h$.
What is the probability of this occurring?

Consider some pair $\{\pi(x_j),\pi(y_j)\}$.
Since $H[K]$ is a random graph on $k$ vertices, their expected agreement set in $H[K]$ is a random variable $X$ with distribution $Bin(k-2,1/2)$.
The probability that it is larger than $0.7h$ is bounded by the standard Chernoff bound (see appendix in \cite{AS-2004}):
$$
\Pr[X > 0.7h]  \le \Pr[X > 0.7(k-2)] < e^{-2 \cdot 0.2^2(k-2)} < e^{-0.06h}
$$
where we have used our assumption that $h \ge k \ge 0.8h$.
Since agreement sets of distinct pairs among $\{\pi(x_1),\pi(y_1)\},\ldots,\{\pi(x_\ell),\pi(y_\ell)\}$ are pairwise independent,
we have from the last inequality and from (\ref{e:l}) that the probability that $\pi$ is an isomorphism between $B_M$ and $H[K]$ is at most 
$$
e^{-0.06h\ell}   < e^{-0.002h^2}\;.
$$
Now, to complete the proof of the lemma, we need to show that {\em no matter} which set of at most $\lfloor 10^{-5}h^2 \rfloor$ pairs we take
to play the role of $M$, and no matter which bijection $\pi$ we take, we still have that $\pi$ is not an isomorphism between $B_M$
(regardless of any additional changes to part pairs of $B$) and $H[K]$, with high probability.
How many possible sets $M$ can we have? As we are selecting at most $\lfloor 10^{-5}h^2 \rfloor$ pairs from a pool of
at most $\binom{h}{2} < h^2/2$ possible pairs, the number of possible $M$ is smaller than
\begin{equation}\label{e:num-m}
\sum _{i=0}^{\lfloor 10^{-5}h^2 \rfloor} \binom{h^2/2}{i} < 
h^2 \binom{h^2/2}{\lfloor 10^{-5}h^2 \rfloor} < e^{0.00007h^2}
\end{equation}
where we have used the standard Stirling approximation.
It follows from the union bound and the fact that there are less than $h!$ possible $\pi$,
that the probability that $B$ is blowup $10^{-5}h^2$-close to $H[K]$ is at most
$$
h! e^{0.00007h^2} \cdot e^{-0.002h^2} < e^{-0.001h^2}\;.
$$
\qed

\begin{lemma}\label{l:cond-5}
If $h$ is sufficiently large and $H \sim G(h,1/2)$, then the following holds.
Let $K \subset [h]$ and $J \subset [h]$ with $|K|=|J|=\lceil 0.7h \rceil$ and let $\pi$ be a bijection from $J$ to $K$
having at least $0.1h$ non-stationary points.
The probability that $\pi$ is $10^{-5}h^2$-close to an isomorphism between $H[J]$ and $H[K]$ is at most $e^{-0.12h^2}$. 
\end{lemma}
\Proof
Let $k=\lceil 0.7h \rceil$ be the order of $J$ and $K$.
Let $A \subset J$ be a set of $s=\lceil 0.1h \rceil$ non-stationary points
and let $B=\pi(A) \subset K$ be  the set of their images. Then, $\pi$ restricted to $A$ is a bijection from $A$ to $B$ with no
stationary points. By Lemma \ref{l:non-stationary}, we can pick $t=\lceil s/3 \rceil$ pairs
$(a_1,\pi(a_1)),\ldots,(a_t,\pi(a_t))$ such that all $2t$ elements in these pairs are distinct.

As in the proof of the previous lemma, let us fix a set $M$ of $m \le \lfloor 10^{-5}h^2 \rfloor$ pairs of vertices of $J$.
Let $H[J]_M$ be the graph obtained from $H[J]$ after changing $M$ (making pairs in $M$ that are edges into non-edges and making
pairs in $M$ that are non-edges into edges).
Call a vertex of $J$ {\em problematic} if it appears in at least $h/30$ elements of $M$.
Observe that being problematic has nothing to do with the structure of $H[J]$, as it only depends on the choice of $M$ and the choice of subset $J$.
As in (\ref{e:num-prob}), the number of problematic vertices is at most $0.0006h$.

Now, remove from the set of pairs $(a_1,\pi(a_1)),\ldots,(a_t,\pi(a_t))$ all the pairs for which $a_i$ is problematic.
The number of pairs that remain is at least
$$
t-0.0006h = \lceil s/3 \rceil - 0.0006h \ge 0.1h/3 - 0.0006h \ge 0.03h\;.
$$
Without loss of generality, let these pairs be $(a_1,\pi(a_1)),\ldots,(a_\ell,\pi(a_\ell))$ where $\ell = \lceil 0.03h \rceil$.
So, all of these $2\ell$ vertices are distinct, and all the $a_1,\ldots,a_\ell$ are non-problematic.
Consider the set $S = \{a_1,\pi(a_1),\ldots,a_\ell,\pi(a_\ell)\}$. Notice that $|J \setminus S| \ge k - 2\ell$.
Now, let $w \in J \setminus S$ be such that $(a_i,w) \notin M$. Since $a_i$ is non-problematic the number of choices for $w$
is at least
$$
|J \setminus S| - h/30 \ge k-2\ell - h/30 \ge 0.7h -2\lceil 0.03h \rceil - h/30 \ge 0.6h\;.
$$
consider the two pairs
$\{a_i,w\}$ and $\{\pi(a_i),\pi(w)\}$. In order for $\pi$ to be an isomorphism between $H[J]_M$ and $H[K]$
we must have that $\{a_i,w\}$ and $\{\pi(a_i),\pi(w)\}$ agree (both are edges or both are non-edges).
Since agreement occurs with probability $\frac{1}{2}$ and since all the choices of $i$ and $w$
result in independent events, since they correspond to distinct pairs, we obtain that
\begin{equation}\label{e:pi-2}
\Pr[\pi \text{ is an isomorphism between $H[J]_M$ and $H[K]$}] \le 2^{-\ell \cdot 0.6h} \le 2^{-0.18h^2}\;.
\end{equation}
We note that a similar idea to the one in the last inequality has been used in the proof of Theorem 3.1 in \cite{KSV-2002}.
To complete the proof of the lemma, we must again, as in the previous lemma, consider all possible choices for $M$.
As in (\ref{e:num-m}), the number of possible $M$ is smaller than $e^{0.00007h^2}$.
It follows from the union bound that the probability that $\pi$ is $10^{-5}h^2$-close to an isomorphism between $H[J]$ and $H[K]$ is at most
$$
e^{0.00007h^2} \cdot 2^{-0.18h^2} < e^{-0.12h^2}\;.
$$
\qed

Recall that in the proof of Lemma \ref{l:cond-5} we have used the following lemma.
\begin{lemma}\label{l:non-stationary}
Let $f$ be a bijection between sets $A$ and $B$ of size $s$ and with no stationary points.
Then there are at least $t=\lceil s/3 \rceil$ pairs $(a_1,f(a_1)),\ldots,(a_t,f(a_t))$ such that all
$2t$ elements in these pairs are distinct.
\end{lemma}
\Proof
A {\em chain} of $f$ is a maximal length sequence $(x_1,\ldots,x_\ell)$ such that $f(x_i)=x_{i+1}$ for $i=1,\ldots\ell-1$ and for which $x_\ell \notin A$ or $f(x_\ell) = x_1$
(in the latter case, the chain forms a {\em cycle}). Clearly, $A \cup B$ can be partitioned into chains. Furthermore, no chain is a singleton since there are no stationary points.
So, from each chain $(x_1,\ldots,x_\ell)$ we can pick at least $\lfloor \ell/2 \rfloor$ pairs $(x_1,x_2=f(x_1))$, $(x_3,x_4=f(x_3))$ and so on, such that all the elements in
the picked pairs are distinct. Thus, if the partition has $c$ odd length chains, we have picked $(|A \cup B|-c)/2$ pairs.
Now, since $|A \cup B| \ge s$ and since $c \le |A \cup B|/3$ the lemma follows. \qed

\noindent
{\bf Proof of Theorem \ref {t:phe}.}
Assume that $h$ is sufficiently large to satisfy all five lemmata \ref{l:cond-1} until \ref{l:cond-5}.
Consider $H \sim G(h,1/2)$. The probability that it does not satisfy the first condition for ${\cal P}_h$ is at most $1/(5h)$ by Lemma \ref{l:cond-1}.
The probability that it does not satisfy the second condition for ${\cal P}_h$ is at most $1/(5h)$ by Lemma \ref{l:cond-2}.
The probability that it does not satisfy the third condition for ${\cal P}_h$ is at most $1/h$ by Lemma \ref{l:cond-3}.
There are less than $2^h$ options for $S \subset [h]$, at most $2^h$ options for $K \subseteq  [h]$ and less than
$h^h$ options for $B$ to be the blowup of $H[S]$. Hence, the probability that $H$ does not satisfy the fourth condition for ${\cal P}_h$
is, by Lemma \ref{l:cond-4} at most
$$
2^h \cdot 2^h \cdot h^h \cdot e^{-0.001h^2} < \frac{1}{5h}
$$
(where here we assume that $h$ is also sufficiently large to satisfy the last inequality).
There are less than $2^h$ options for each of $J$ and $K$ in the fifth condition of ${\cal P}_h$, and less than $h!$ possible bijections $\pi$.
Hence, the probability that $H$ does not satisfy the fifth condition for ${\cal P}_h$
is, by Lemma \ref{l:cond-5} at most
$$
2^h \cdot 2^h \cdot h! \cdot e^{-0.12 h^2} < \frac{1}{5h}\;.
$$
Hence the probability that $H \notin {\cal P}_h$ is at most
$$
4 \cdot \frac{1}{5h} + \frac{1}{h} < \frac{2}{h}\;.
$$
\qed

\section{Proof of the main result}

We will prove a stronger theorem which immediately implies Theorem \ref{t:main} by just using $\epsilon = \frac{1}{2}$ in the following statement and by
redefining  ${\cal P}_h = \emptyset$ for $h \le h_0(\frac{1}{2})$.
\begin{theorem}\label{t:super}
For every $\epsilon > 0$ there exists $h_0=h_0(\epsilon)$ such that the following hold for all $h > h_0$
and for all $n \le 2^{h^{1-\epsilon}}$.
\begin{enumerate}
\item
If $H \in {\cal P}_h$, then $i_H(n)=g(n,h)$. Furthermore, the family of extremal graphs is precisely $H^*(n)$.
\item
if $G$ has $n$ vertices and $i_H(G) \ge f(n,h)$ then $G$ must be obtained from a balanced blowup of $H$ by adding some edges inside the blowup parts.
\end{enumerate}
\end{theorem}
\Proof
The fact that $i_H(n) \ge g(n,h)$ was already shown in Section 2. Thus, we must prove that $i_H(n) \le g(n,h)$
and that the only graphs which attain this bound are the nested balanced blowups of $H$.
Throughout the proof we shall assume that $h$ is sufficiently large as a function of $\epsilon$ to satisfy various inequalities,
thereby establishing the existence of $h_0$. The proof will be by induction on $n$ where the base case $n \le h$ is trivial since
$g(h,h)=1$ and since $g(n,h)=0$ for all $n < h$.

Let $n \le 2^{h^{1-\epsilon}}$ and assume that we have already established the result for all values smaller than $n$.
Let $G$ be a graph with $n$ vertices. Let $\H$ denote the set of all induced copies of $H$ in $G$.
In other words, considering $G$ and $H$ as labeled graphs,
each $F \in \H$ is identified with an injective function from $V(H)$ to $V(G)$ which preserves adjacencies
and non-adjacencies (thus, we view $F$ as such a function). We must therefore prove that $|\H| \le g(n,h)$.

Assume that $V(H)=[h]$. Since $H \in {\cal P}_H$, it has a distinguishing set $Q \subset [h]$ with
$|Q|=q \le 3\log h$. Recall that $Q$ has the following property.
For any two distinct vertices $u,v \in [h] \setminus Q$, there exists $w \in Q$
such that exactly one of $u$ or $v$ is adjacent to $w$. For the rest of this proof we fix $Q$.

\begin{definition}{\bf [Role consistent]}
Let $F_1,F_2 \in \H$ be two induced copies of $H$ in $G$. We say that $F_1$ and $F_2$ are
{\em role-consistent} if for all $x \in Im(F_1) \cap Im(F_2)$ we have $F_1^{-1}(x) = F_2^{-1}(x)$.
A subset $\F \subseteq \H$ is {\em role-consistent} if any pair of
elements of $\F$ are role-consistent.
\end{definition}

A role-consistent subset $\F$ induces a partition $P(\F) = \{R,P_1,\ldots,P_h\}$ of $V(G)$ as follows.
For each $v \in V(G)$, if $v$ is in no element of $\F$ (namely, in no image of any $F \in \F$),
then assign $v$ to $R$. Otherwise, let $v \in Im(F)$ for some $F \in \F$.
Then, place $v$ in $P_j$ where $F^{-1}(v)=j$. Since $\F$ is role-consistent, this assignment is consistent regardless of
the choice of $F \in \F$ for which $v \in Im(F)$.
Stated otherwise, every copy of $H$ in $\F$ has the property that it contains precisely one vertex from each of
$P_1,\ldots,P_h$ where the vertex of the copy that plays the role of $i$ belongs to $P_i$.

\begin{definition}{\bf [Role partition; redundant part]}
We call $P(\F)$ the {\em role partition} of $\F$ and call $R$ the {\em redundant part} of the role partition.
\end{definition}

It immediately follows that if $\F$ is role-consistent and  $P(\F) = \{R,P_1,\ldots,P_h\}$ is its role partition, then
\begin{equation}\label{e:fnh}
|\F| \le \Pi_{i=1}^h |P_i| \le f(n-|R|,h) \le f(n,h)\;.
\end{equation}

\begin{definition}{\bf [$Q$-partition]}
Let $\I$ be the set of all injective functions from $Q$ to $V(G)$. For $I \in \I$ let $\H_I$ be the set of all
elements $F \in \H$ which are consistent with $I$. Namely, $F(x)=I(x)$ for all $x \in Q$.
Clearly $\{\H_I ~|~ I \in \I\}$ is a partition of $\H$ which we call the {\em $Q$-partition}.
\end{definition}
Observe also that since $|Q|=q$ we have $|\I| < n^q$. 

Our high level approach is to prove that a delicate modification of the $Q$-partition has the property that most
of its elements are small, and that the remaining (not so small) elements of this partition can be absorbed into one large part
that is still role-consistent. We first need the following lemma.

\begin{lemma}\label{l:main-1}
Each element of the $Q$-partition is role-consistent. Namely, for each $I \in \I$, the set $\H_I$ is role-consistent.
\end{lemma}
\Proof
Let $F_1,F_2 \in \H_I$. Let $x \in Im(F_1) \cap Im(F_2)$.
Suppose $F_1^{-1}(x)=i$ and $F_2^{-1}(x)=j$. Now, if $i \in Q$ or $j \in Q$, then by the definition of
$\F_I$ we must have $i=j$. Otherwise, if $i \notin Q$ and $j \notin Q$, then for every
vertex $w \in Q$, either both $i$ and $j$ are adjacent to $w$ or both are non-adjacent.
Thus, by the property of $Q$, $i$ and $j$ must be the same vertex of $H$, namely $i=j$.
\qed

By Lemma \ref{l:main-1}, by (\ref{e:fnh}) and by the fact that $|\I| < n^q$ we have the following crude upper bound.
\begin{equation}\label{e:crude}
|\H| \le f(n,h)n^q \le f(n,h) n^{3 \log h}\;.
\end{equation}
Notice that this bound applies to all $n \ge h$. Recall however, that the upper bound obtained in (\ref{e:crude}) is still far from what we require
as we would like to essentially eliminate the $n^{3 \log h}$ factor.

Consider the following process performed on a role-consistent set $\F$ which we call {\em core production}.
Start by defining  $X_0=\cal F$. As long as $X_i$ has the property that there is a vertex $v$ of $G$
that appears in at most $f(n,h)n^{-3 -3 \log h}$ elements of $X_i$ ($X_i$ is a subset of $\F$ hence elements of $X_i$ are induced copies of $H$ in $G$), then remove from $X_i$ all elements containing $v$
and denote the remaining set by $X_{i+1}$. When the process ends at some stage $t$ we either have $X_t=\emptyset$ or else $X_t$ is a subset of
$\F$ having the property that each vertex appearing in any element of $X_t$ appears in at least $f(n,h)n^{-3 -3 \log h}$ elements of $X_t$.
Set $X_t = \F^*$.
\begin{definition}{\bf [core; leftover; large]}
Let $\F$ be role-consistent. We call $\F^*$ a {\em core} of $\F$ and call $\F^\# = \F \setminus \F^*$ a {\em leftover} of $\F$.
We say that a nonempty $\F$ is {\em large} if $\F^*=\F$.
\end{definition}
Notice that since $\F^*$ is a subset of $\F$, it is also role-consistent. Also observe that $(\F^*)^*=\F^*$ so every nonempty core is large.
Finally, observe that if $\F$ is large, then every vertex of $G$ that appears in an element of $\F$ appears in at least $f(n,h)n^{-3 -3 \log h}$ other elements,
so in particular $|\F| \ge f(n,h)n^{-3 -3 \log h}$.

\begin{lemma}\label{l:main-2}
Suppose we performed core production on each element of the $Q$-partition, thereby obtaining cores $\{\H_I^*~|~ I \in \I\}$ and leftovers $\{\H_I^\#~|~ I \in \I\}$.
Let $\H^\# = \cup_{I \in \I} \H_I^\#$.
Then,
$$
|\H^\#| \le \frac{f(n,h)}{n^2}\;.
$$
\end{lemma}
\Proof
In the core production of a role-consistent $\F$ we add at each step at most $f(n,h)n^{-3 -3 \log h}$ elements to the leftover set,
where all the added elements of the step contain some vertex $v \in V(G)$.
Since there are $n$ vertices, the final leftover of $\F$ satisfies
$$
|\F^\#| \le n \cdot f(n,h)n^{-3 -3 \log h} = f(n,h)n^{-2 -3 \log h}\;.
$$
Since the number of elements $I \in \I$ is less than $n^q \le n^{3\log h}$, we have
$$
|\H^\#| \le \sum_{I \in \I} |\H_I^\#| \le n^{3\log h}f(n,h)n^{-2 -3 \log h}=\frac{f(n,h)}{n^2}\;.
$$
\qed

Our main lemma (whose proof is deferred to Section \ref{sec:main-lemma}) is the following.
\begin{lemma}\label{l:main-3}
Suppose that $\F$ and $\E$ are two large sets.
Then their union is role-consistent.
\end{lemma}

Once we establish Lemma \ref{l:main-3} we can complete the proof of Theorem \ref{t:super} as follows.
Consider the set of cores $\{\H_I^*~|~ I \in \I\}$ and leftovers $\{\H_I^\#~|~ I \in \I\}$.
Recall that if $\H_I^* \neq \emptyset$, then $\H_I^*$ is large. Hence, by Lemma  \ref{l:main-3} used repeatedly,
$$
\H^* = \cup_{I \in {\cal I}} \H_I^*
$$
is role-consistent (observe that if $\F$ and $\E$ are large and their union is role-consistent then this union is large as well).
There are now two cases to consider:

\noindent
{\bf Case 1:} $|\H^*| < f(n,h)$.
Consider the role-partition $P(\H^*)=\{R,P_1,\ldots,P_h\}$.
Now, $|\H^*| < f(n,h)$ could happen for three reasons.
(i) Either $R \neq \emptyset$. (ii) Else, the partition is not equitable.
(iii) Else, the partition is equitable, but there is some pair of vertices $u,v$ in $G$ such that $u \in P_i$, $v \in P_j$, $i \neq j$,
and the following occurs. either $uv \in E(G)$ but $ij \notin E(H)$ or $uv \notin E(G)$ but $ij \in E(H)$.
Equivalently, no element of $\H^*$ contains both $u$ and $v$.

Suppose (i) occurs. In this case, $|\H^*| \le f(n-1,h)$. But then by Lemma \ref{l:simple} and Lemma \ref{l:main-2} we obtain:
$$
|\H| \le |\H^*| + |\H^\#| \le f(n-1,h)+\frac{f(n,h)}{n^2} = \left(1-\frac{1}{\lceil n/h \rceil}\right)f(n,h)+\frac{f(n,h)}{n^2} < f(n,h)\;.
$$

Suppose (ii) occurs. In this case, $|\H^*|$ is at most $\Pi_{i=1}^h|P_i|$ and the partition is not equitable, hence by Lemma \ref{l:simple} and Lemma \ref{l:main-2} we obtain:
$$
|\H| \le |\H^*| + |\H^\#| \le f(n,h)\left(1-\frac{1}{\lceil n/h \rceil^2}\right)+\frac{f(n,h)}{n^2} < f(n,h)\;.
$$

Suppose (iii) occurs. In this case, no element of $\H^*$ contains both $u$ and $v$ and the partition is equitable.
Thus,
$$
|\H^*| \le f(n,h)\left(1-\frac{1}{|P_i||P_j|}\right) \le f(n,h)\left(1-\frac{1}{\lceil n/h \rceil^2}\right)\;.
$$
So, again by Lemma \ref{l:main-2} we obtain:
$$
|\H| \le |\H^*| + |\H^\#| \le f(n,h)\left(1-\frac{1}{\lceil n/h \rceil^2}\right) + \frac{f(n,h)}{n^2} < f(n,h)\;.
$$

\noindent
{\bf Case 2:} $|\H^*| = f(n,h)$. First observe that the structure of $\H^*$ is unique. Clearly we must have $P(\H^*)=\{\emptyset,P_1,\ldots,P_h\}$
and the partition is equitable. Furthermore, for each choice of $h$ vertices, one from each part (i.e. a transversal), we obtain an element of $\H^*$,
as this is the only way to get $f(n,h)$ elements in $\H^*$. Observe also that this means that the bipartite graph between $P_i$ and $P_j$ is either empty
(if $ij \notin E(H)$) or complete (if $ij \in E(H)$).
Notice that as a bonus we have now proved the second part of Theorem \ref{t:super}, since in the previous case where  $|\H^*| < f(n,h)$ we have already shown that
$|\H| < f(n,h)$.

We claim that every element of $\H^\#$ is entirely contained in some $P_i$.
Assume otherwise, and let $X \in H^\#$ have its vertices in more than one part.
Observe that it cannot have each vertex in a distinct part, as in this case we are repeating a copy of $H$ which is already in $\H^*$, while $\H^* \cap H^\# = \emptyset$.
So, $X$ has vertices in more than one part but not in all parts.
Consider the smallest non-singleton part, say $P_i$. Thus, $|P_i \cap V(X)| \ge 2$.
If $|P_i \cap V(X)| < 0.4h$, then any two vertices from $P_i \cap V(X)$ have all the vertices in $V(X) \setminus P_i$, namely more than $0.6h$ vertices, in their
agreement set. This contradicts the assumption that $H \in {\cal P}_H$.
If $|P_i \cap V(x)| > 0.6h$ then any vertex in $V(X) \setminus P$ has either all the vertices of $P_i \cap V(x)$ as its neighbors, or as its non-neighbors.
But this contradicts the assumption that for $H \in {\cal P}_H$, the degree of any vertex is between $0.4h$ and $0.6h$.
Finally, if $0.4h \le |P_i \cap V(X)| \le 0.6h$ then any three vertices from $P_i \cap V(X)$ have all the vertices in $V(X) \setminus P_i$, namely more than $0.4h$ vertices, in their
agreement set, contradicting the assumption that $H \in {\cal P}_H$.

Having proved that every element of $\H^\#$ is entirely contained in some $P_i$, we can now use the induction hypotheses to obtain that
$$
|\H^\#| \le \sum_{i=1}^h g(|P_i|,h).
$$
Therefore,
$$
|\H| \le |\H^*| + |\H^\#| \le f(n,h) + \sum_{i=1}^h g(|P_i|,h) = g(n,h)\;.
$$
Finally, notice that the only way to obtain an equality is to
have each $P_i$ induce a graph with the maximum number of induced copies of $H$, namely, by induction, with $g(|P_i|,h)$ induced copies of $H$.
Thus, each $P_i$ induces an element of the nested balanced blowup on $|P_i|$ vertices. Hence, by the definition of $H^*(n)$, $G$ is
a nested balanced blowup of $H$ with $n$ vertices, namely $G \in H^*(n)$.
\qed

\section{Properties of large sets}\label{sec:properties}

In order to prove Lemma \ref{l:main-3}, we first need to establish some properties that hold for every large set.
Let $\F$ be a large set and let $P(\F)=\{R,P_1,\ldots,P_h\}$ be the role partition of $\F$.
Our first lemma quantifies the fact that in a large set, a small fraction of parts cannot contain too many elements.

\begin{lemma}\label{l:some-parts}
Let $J \subset [h]$ with $|J| = \lfloor 0.00001h \rfloor$. Then, $\sum_{j \in J} |P_j| \le 0.0004n$.
\end{lemma}
\Proof
Consider the following process. Let $n \ge h$. Start with an equitable partition of $n$ elements into $h$ parts.
The product of the sizes of the parts is thus $f(n,h)$. Let $0 < k < h$. Designate the largest $k$ parts of the equitable partition
(so we prefer take parts of size $\lceil n/h \rceil$ to the designated parts as long as there are enough, and if there are less than $k$ such parts,
then the remaining designated parts are of size $\lfloor n/h \rfloor$). Now, suppose $t \ge k \lceil n/h \rceil$.
Note that the designated parts contain altogether at most $t$ elements.
Move elements from the non-designated parts to the designated parts so that after each move, the designated parts form an equitable partition
of $k$ parts and the non-designated parts from an equitable partition into $h-k$ parts. Stop moving after the designated parts contain precisely $t$ elements.
The number of moves that we have made is therefore at least
$$
t- \lceil n/h \rceil \cdot k\;.
$$
What happened to the product of the sizes of the parts after each move?
Suppose we have moved an element from a non-designated part whose current size is $a$ to a designated part whose current size is $b$,
hence $b \ge a$.
So the ratio between the product of the part sizes after the move and the product before the move is precisely
$$
\frac{(a-1)(b+1)}{ab} = \left(1-\frac{b-a+1}{ab} \right) < 1-\frac{1}{a}+\frac{1}{b}\;.
$$
Notice that at any point, we have $a \le \lceil n/h \rceil$. Also, after half of the moves are made, each element in a
designated set contains at least $\lfloor t/2k \rfloor$ elements.
Hence for the last half of the moves, we have $b \ge \lfloor t/2k \rfloor$.
In these moves the aforementioned ratio is at most
$$
1-\frac{1}{\lceil n/h \rceil}+ \frac{1}{\lfloor t/2k \rfloor}\;.
$$
So after all moves are completed, the product of the sizes of the parts is now at most
\begin{equation}\label{e:process}
f(n,h) \left(1-\frac{1}{\lceil n/h \rceil}+ \frac{1}{\lfloor t/2k \rfloor}\right)^{t/2-\frac{\lceil n/h \rceil \cdot k}{2}}\;.
\end{equation}

Let us now return to the statement of the lemma.
We will consider the product of the sizes of the parts $\Pi_{i=1}^h |P_i|$.
Let $k=\lfloor 0.00001h \rfloor$ and let $t > 0.0004n$.
Suppose the lemma does not hold. Then we have $k$ parts that together contain $t$ elements.
By convexity, the maximum of the product of the sizes of the parts subject to this condition is obtained when $R = \emptyset$ (no redundant
part), and when the $k$ parts form an equitable partition of $t$ into $k$ parts, and the remaining $h-k$ parts form an equitable partition
of $n-t$ into $h-k$ parts. By (\ref{e:process}) we have that
\begin{eqnarray*}
\Pi_{i=1}^h |P_i| & \le & f(n,h) \left(1-\frac{1}{\lceil n/h \rceil}+ \frac{1}{\lfloor t/2k \rfloor}\right)^{t/2-\frac{\lceil n/h \rceil \cdot k}{2}}\\
& \le & f(n,h) \left(1-\frac{1}{\lceil n/h \rceil}+ \frac{1}{\lfloor 5n/h \rfloor}\right)^{0.0002n-\frac{\lceil n/h \rceil \cdot \lfloor 0.00001h \rfloor}{2}}\\
& \le & f(n,h) \left(1-\frac{h}{4n}\right)^{0.0001n}\\
& \le & f(n,h) e^{-0.000025h}\\
& \le & f(n,h) 2^{-h^{1-\epsilon}4\log h}\\
&  <  & f(n,h) n^{-4 \log h}\;.
\end{eqnarray*}
Here we have used $t/(2k) \ge 5n/h$ and $\frac{1}{\lceil n/h \rceil}- \frac{1}{\lfloor 5n/h \rfloor} \ge h/(4n)$
while recalling also our assumptions that $h$ is sufficiently large as a function of $\epsilon$ and that $n \le 2^{h^{1-\epsilon}}$.
Now, since $|\F| \le \Pi_{i=1}^h |P_i|$ we have reached a contradiction to the assumption that $\F$ is large,
as large sets have at least $f(n,h) n^{-3-3 \log h}$ elements.
\qed

Call a part $P_i$ {\em small} if it contains at most $\lfloor n/(10h) \rfloor$ vertices.
\begin{lemma}\label{l:small-parts}
The number of small parts is less than $\lfloor 0.0002h \rfloor$.
\end{lemma}
\Proof
The proof is similar to the previous one, only now we wish to focus on moving elements away from a set of designated parts
so the bounds we obtain on the product of the set sizes fits better to this scenario.
Consider the following process. Let $n \ge h$. Start with an equitable partition of $n$ elements into $h$ parts.
The product of the sizes of the parts is thus $f(n,h)$. Let $0 < k < h$. Designate the smallest $k$ parts of the equitable partition.
Now, suppose $2t \le k \lfloor n/h \rfloor$. Note that the designated parts contain altogether at least $2t$ elements.
Move elements from designated parts to the non-designated parts so that after each move, the designated parts form an equitable partition
of $k$ parts and the non-designated parts from an equitable partition into $h-k$ parts.
Stop moving after the designated parts contain precisely $t$ elements.

What happened to the product of the sizes of the parts after each move?
Suppose we have moved an element from a designated part whose current size is $a$ to a non-designated part whose current size is $b$,
hence $b \ge a$.
So the ratio between the product of the part sizes after the move and the product before the move is smaller than $1-\frac{1}{a}+\frac{1}{b}$.
Notice that at any point, we have $b \ge \lfloor n/h \rfloor$. Consider the last $t$ moves. In each such move, the total number of elements in the designated sets
prior to the move is at most $2t$. Hence, for each such move we have $a \le \lceil 2t/k \rceil$.
So in these moves  the aforementioned ratio is at most
$$
1-\frac{1}{\lceil 2t/k \rceil}+ \frac{1}{\lfloor n/h \rfloor}\;.
$$
So after all moves are completed, the product of the sizes of the parts is now at most
\begin{equation}\label{e:process-2}
f(n,h) \left(1-\frac{1}{\lceil 2t/k \rceil}+ \frac{1}{\lfloor n/h \rfloor}\right)^t\;.
\end{equation}

Let us now return to the statement of the lemma.
We will consider the product of the sizes of the parts $\Pi_{i=1}^h |P_i|$.
First, observe that we may assume that $n \ge 10h$ as otherwise having
a small part means having an empty part, meaning that $\F = \emptyset$ which is impossible as we assume that
$\F$ is large.

Let $k=\lfloor 0.0002h \rfloor$.
Suppose the lemma does not hold. Then we have $k$ small parts that together contain $t$ elements.
By convexity, the maximum of the product of the sizes of the parts subject to this condition is obtained when $R = \emptyset$ (no redundant
part), when the $k$ parts form an equitable partition of $t$ into $k$ parts, the remaining $h-k$ parts form an equitable partition
of $n-t$ into $h-k$ parts, and $t$ is maximized namely $t=k\lfloor n/10 h \rfloor$. By (\ref{e:process-2}) we have that
\begin{eqnarray*}
\Pi_{i=1}^h |P_i| & \le & f(n,h) \left(1-\frac{1}{\lceil 2 \lfloor n/10 h \rfloor \rceil}+ \frac{1}{\lfloor n/h \rfloor}\right)^{\lfloor 0.0002h \rfloor\lfloor n/10 h \rfloor}\\
& \le & f(n,h) \left(1-\frac{h}{n}\right)^{10^{-5}n}\\
& \le & f(n,h) e^{-10^{-5}h}\\
& \le & f(n,h) 2^{-h^{1-\epsilon}4\log h}\\
&  <  & f(n,h) n^{-4 \log h}\;.
\end{eqnarray*}
Here we have used $\frac{1}{\lceil 2 \lfloor n/10 h \rfloor \rceil}- \frac{1}{\lfloor n/h \rfloor} \ge h/n$
while recalling also our assumptions that $h$ is sufficiently large as a function of $\epsilon$ and that $n \le 2^{h^{1-\epsilon}}$.
Now, since $|\F| \le \Pi_{i=1}^h |P_i|$ we have reached a contradiction to the assumption that $\F$ is large,
as large sets have at least $f(n,h) n^{-3-3 \log h}$ elements.
\qed

Let $u \in V(G) \setminus R$ and let ${\hat \F}_u$ be the set of elements of $\F$ that contain $u$.
Let $B(u,\F) \subset V(G)$ be the set of vertices of $G$ that do not appear in any element of ${\hat \F}_u$.
Clearly, $R \subseteq B(u,\F)$ as vertices in the redundant part do not appear in any element of $\F$.
Since $\F$ is large, we have that any vertex of $G$ which is not in $R$, appears in at least $f(n,h)n^{-3 - 3\log h}$ elements of $\F$.
Thus, $|{\hat \F}_u| \ge f(n,h)n^{-3 - 3\log h}$.
Let $\F_u$ be the set of elements of $\F$ which do not contain vertices from $B(u,\F)$.
Notice that the role partition of $\F_u$ is precisely $P(\F_u)=\{R',P_1',\ldots,P_h'\}$ where $R' = B(u,\F)$ and $P_i' = P_i \setminus B(u,\F)$.
Also, clearly, ${\hat \F}_u \subseteq \F_u$ but the inclusion may be proper as there may be elements in $\F$ that do not contain $u$ nor any element from $B(u,\F)$.
Given this inclusion, we have
$$
|\F_u| \ge |{\hat \F}_u| \ge f(n,h)n^{-3 - 3\log h}\;.
$$
Our next lemma bounds the order of $B(u,\F)$.

\begin{lemma}\label{l:num-bad}
$|B(u,\F)| \le 10^{-5}n$.
\end{lemma}
\Proof
Suppose $|B(u,\F)|=k$.
Since each element of ${\hat \F}_u$ contains no vertex of $B(u,\F)$, we have that $|{\hat \F}_u| \le f(n-k,h)$.
Using Lemma \ref{l:simple} we have
$$
|{\hat \F}_u|\le f(n-k,h) \le f(n,h) \left(1-\frac{1}{\lceil n/h \rceil} \right)^k\;.
$$
We must show that $k \le 10^{-5}n$. Assume otherwise, then
$$
\frac{|{\hat \F}_u|}{f(n,h)}  \le \left(1-\frac{1}{\lceil n/h \rceil} \right)^k \le \left(1-\frac{h}{2n} \right)^{10^{-5}n} < e^{-10^{-6}h} < 2^{-h^{1-\epsilon}4\log h} < n^{-4 \log h}
$$
contradicting the fact that $|{\hat \F}_u| \ge f(n,h)n^{-3 - 3\log h}$.
\qed

Call a part $P'_i$ {\em very small} if it contains at most $n/(15h)$ vertices.
\begin{lemma}\label{l:very-small}
The number of very small parts is at most $0.0005h$.
\end{lemma}
\Proof
By Lemma \ref{l:small-parts}, the number of small parts $P_i$ is less than $\lfloor 0.0002h \rfloor$. How many parts $P_i$ that were not small
turned out to be very small parts in $P_i' = P_i \setminus B(u,F)$. For this to happen, $P_i$ must lose at least $n/(10h) - n/(15h) = n/(30h)$
of its vertices because they joined $B(u,F)$. But since by Lemma \ref{l:num-bad} $|B(u,\F)| \le 10^{-5}n$, we have that the number of
such parts is at most $(n/100000)/(n/(30h)) = 0.0003h$. Hence, overall there are at most $0.0003h+\lfloor 0.0002h \rfloor \le 0.0005h$ very small parts.
\qed

Let $u \in P_i'$ and $v \in P_j'$ with $i \neq j$. We say that the pair is {\em inconsistent} if either
$ij \in E(H)$ and $uv \notin E(G)$ or $ij \notin E(H)$ and $uv \in E(G)$. Otherwise, the pair is {\em consistent}.
Let $z > 0$ be a real parameter. We define the $z$-consistency graph $Z_{\F_u}(z)$ as follows.
The set of vertices of $Z_{\F_u}(z)$ is $[h]$ (namely, it has the same set of vertices as $H$).
Two vertices $i$ and $j$ are adjacent in $Z_{\F_u}(z)$ if and only if the number of consistent pairs with one endpoint in $P_i'$ and one endpoint in $P_j'$ is
more than $(1-z)|P'_i||P'_j|$.

Our next lemma shows that for relatively small $z$, $Z_{\F_u}(z)$ has large cliques.
\begin{lemma}\label{l:zfz}
$Z_{\F_u}(h^{-\epsilon/2})$ has a clique of size at least $0.99999h+1$.
\end{lemma}
We will actually prove that the complement of $Z_{\F_u}(h^{-\epsilon/2})$ does not have a matching of size at least $0.000005h-1$, as this implies that
$Z_{\F_u}(h^{-\epsilon/2})$ has a clique of size at least $h-2\cdot 0.000005h+2=0.99999h+2$. Denote this complement by $C$.

Suppose $C$ has a matching of size $k=\lceil 0.000005h \rceil-1$. Denote it by $(x_1,y_1), \ldots, (x_k,y_k)$.
Any copy of $H$ in $\F_u$ contains precisely one vertex from each of the parts $P_1',\ldots,P_h'$.
By the definition of $Z_{\F_u}(z)$, each such copy must avoid at least $z|P'_{x_i}||P'_{y_i}|$ of the pairs of vertices from $P'_{x_i} \times P'_{y_i}$,
as it must avoid inconsistent pairs and this holds for all $i=1,\ldots,k$. Hence, for $z=h^{-\epsilon/2}$ we have
\begin{eqnarray*}
|\F_u| & \le & (1-z)^k \Pi_{i=1}^h |P'_i|\\
& \le & \left( 1-\frac{1}{h^{\epsilon/2}}\right)^{0.000005 h-1} \Pi_{i=1}^h |P'_i|\\
& \le & e^{-0.000004h^{1-\epsilon/2}} f(n,h) \\
& \le & f(n,h)2^{-h^{1-\epsilon}4\log h} \\
& \le & f(n,h) n^{-4\log h}
\end{eqnarray*}
contradicting the fact that $|\F_u| \ge f(n,h)n^{-3 - 3\log h}$.
\qed

As we will be interested in the case $z=h^{-\epsilon/2}$, we will denote for simplicity $Z_{\F_u}=Z_{\F_u}(h^{-\epsilon/2})$.
As Lemma \ref{l:zfz} suggests, we have a clique of order $\lceil 0.99999h \rceil+1$ in $Z_{\F_u}$ so let $K^* \subset [h]$ denote such a clique.
By Lemma \ref{l:very-small}, there are at most $0.0005h$ very small parts. Hence, there is $K(u,\F) \subseteq K^*$ with
$|K(u,\F)| \ge (\lceil 0.99999h\rceil +1)- 0.0005h -1\ge 0.999h$ such that $K(u,\F)$ is a clique in $Z_{\F_u}$ and for any $j \in K(u,\F)$ we have that $|P'_j| \ge n/(15h)$.
Furthermore, if $u \in P_\ell$ then $\ell \notin K(u,\F)$.

How many vertices of $G$ appear in parts $P_i'$ for which $i \in K(u,\F)$ (namely, in the clique parts that are not very small, and not containing the part
of $u$).
We claim that these are most of the vertices. Denote this set by $W(u,\F) = \cup_{i \in K(u,\F)} P'_i$.
\begin{lemma}\label{l:num-good}
$|W(u,\F)| \ge 0.9995 n$\;.
\end{lemma}
\Proof
A vertex $v \notin W(u,\F)$ has to satisfy one of the following conditions. Either $v \in B(u,\F)$, else $v \in P_i'$ but $i \notin K^*$,
else $v \in P_i'$, $i \in K^*$ but $P'_i$ is very small, else $v$ is in the same part as $u$.
Summing the sizes of these forbidden sets we obtain by
Lemmas \ref{l:some-parts}, \ref{l:num-bad}, \ref{l:very-small} and the fact that
$|[h]\setminus K^*| + 1 \le \lfloor 0.00001h \rfloor$ that it is at most
$$
\frac{n}{100000}+0.0004n+0.0005h\cdot (n/(15h)) \le 0.0005 n\;.
$$
\qed

Let $C^*(u,\F)$ denote the set of inconsistent pairs $x,y$ with $x \in P_i'$, $y \in P_j'$ and both $i,j \in K(u,\F)$.
\begin{lemma}\label{l:num-inconsistent}
$|C^*(u,\F)| \le n^2 h^{-\epsilon/2}$.
\end{lemma}
\Proof
Clearly, by the definition of $K(u,\F)$, the number of inconsistent pairs with one endpoint in $P_i'$ and the other in $P_j'$
is at most $|P'_i||P'_j|h^{-\epsilon/2}$ when $i,j \in K(u,\F)$ (in fact, this also holds in the larger set $K^*$).
Hence,
\begin{eqnarray*}
|C^*(u,\F)| & \le & h^{-\epsilon/2} \sum_{i,j \in K(u,\F), i \neq j} |P'_i||P'_j|\\
& \le & h^{-\epsilon/2} \sum_{i=1}^h \sum_{j=i+1}^h |P_i||P_j|\\
& \le & h^{-\epsilon/2} h^2 \frac{n^2}{h^2} \le n^2 h^{-\epsilon/2}\;.
\end{eqnarray*}
\qed

\section{Proof of Lemma \ref{l:main-3}}\label{sec:main-lemma}

We now return to the proof of Lemma \ref{l:main-3}.
Recall that we want to prove that if $\F$ and $\E$ are two large sets, then their union is role-consistent.
Let $P(\F)=\{R,P_1,\ldots,P_h\}$ and $P(\E)=\{T,S_1,\ldots,S_h\}$ be the role partitions of $\F$ and $\E$.
We must therefore prove that if $u \in P_i$ and $u \in S_j$, then $i=j$.

Assume otherwise, then there exists some $u \in V(G)$ such that $u \in P_i$ and $u \in S_j$ and $i \neq j$.
Using the aforementioned notations we focus on the objects
$B(u,\F)$, $\F_u$, $P(\F_u)=\{R',P_1',\ldots,P_h'\}$, $Z_{\F_u}$, $K(u,\F)$, $W(u,\F)$, $C^*(u,\F)$
and correspondingly 
$B(u,\E)$, $\E_u$, $P(\E_u)=\{T',S_1',\ldots,S_h'\}$, $Z_{\E_u}$, $K(u,\E)$, $W(u,\E)$, $C^*(u,\E)$.
Also recall that $R' = B(u,\F)$ and $T'=B(u,\E)$.

Call an index $i \in K(u,\E)$ {\em problematic for $\F$} if more than $|S'_i|/10$ of the vertices of $S'_i$ do not belong to $W(u,\F)$.
Notice that unlike the definitions in the previous section which only depend on a given large set and a vertex of its partition,
the definition of being problematic involves two large sets.
Informally, being non-problematic means that $90$ percent of the vertices of the part appear in ``good sets''
of the other partition. How many problematic parts are there?
\begin{lemma}\label{l:num-problematic}
The number of vertices of $K(u,\E)$ that are problematic for $\F$ is at most $0.075h$.
\end{lemma}
\Proof
Suppose $i \in K(u,\E)$. Then, in particular $S'_i$ is not very small so it contains at least $n/(15h)$ vertices.
If it were problematic for $\F$, then it would have at least $0.1 n/(15h) = n/(150h)$ vertices not in $W(u,\F)$.
But the number of vertices not in $W(u,\F)$ is at most $0.0005n$ by Lemma \ref{l:num-good}.
So, the overall number of vertices of $K(u,\E)$ that are problematic for $\F$ is at most $0.0005n/(n/(150h)) \le 0.075h$.
\qed

Let $J$ denote the set of vertices of $K(u,\E)$ that are non-problematic for $\F$.
Recall that $|K(u,\E)| \ge 0.999h$. Hence, by Lemma \ref{l:num-problematic}, $|J| \ge (0.999-0.075)h \ge 0.9h$.

Let $X \subset V(G)$ be a set of $|J|$ vertices obtained by {\em randomly} selecting precisely one vertex from each part $S'_j$ for $j \in J$.
All $|J|$ choices are performed independently.
Thus, $G[X]$ is a subgraph of $G$ on $|J|$ vertices. We would like to first see how close $G[X]$ is to a subgraph of $H$ on $|J|$ vertices.
\begin{lemma}\label{l:far-x}
With probability larger than $3/4$, $G[X]$ is $4h^{2-\epsilon/2}$-close to $H[J]$.
\end{lemma}
\Proof
Denote $X=\{x_j~|~j \in J\}$. So $x_j$ is the vertex selected at random from $S'_j$.
Consider a pair of distinct indices $i,j$ such that $i,j \in J$. What is the probability that $x_i$ and $x_j$ are consistent?
First notice that since $i,j \in J$, then $|S'_i|,|S'_j| \ge n/(15h)$ as they are not very small.
Next, notice that $ij$ is a clique edge in $Z_{\E_u}$, and hence the number of pairs between $S'_i$ and $S'_j$ that
are consistent is at least $(1-h^{-\epsilon/2})|S'_i||S'_j|$. So the probability that $x_i$ and $x_j$ are inconsistent
is at most
$h^{-\epsilon/2}$. Thus, the expected number of inconsistent pairs in $X$ is at most
$$
\binom{|J|}{2}h^{-\epsilon/2} < h^{2-\epsilon/2}\;.
$$
By Markov's inequality the probability of the number of inconsistent pairs being more than $4h^{2-\epsilon/2}$ is less than $1/4$.
The proof is complete by noticing that changing the inconsistent pairs results in a graph that is isomorphic to $H[J]$ where $x_j$ plays the role of $j$ in $H[J]$.
\qed

\begin{lemma}\label{l:x*}
Let $X^* = X \cap W(u,\F)$. Then with probability larger than  $3/4$, $|X^*| > 0.8h$.
\end{lemma}
\Proof
Recall that $X=\{x_j~|~j \in J\}$ where $x_j$ is the vertex selected at random from $S'_j$.
Since $S'_j$ is non-problematic for $\F$, at least a fraction of $0.9$ of its vertices are in $W(u,\F)$. Hence the probability that
$x_j \in  W(u,\F)$ is at least $0.9$. The expected value of $|X^*|$ is therefore at least $0.9|J| \ge 0.9 \cdot 0.9h = 0.81h$.
As $|X^*|$ is the sum of independent indicator random variables, the probability that it is smaller by a constant factor
than its expectation is exponentially small in $h$. Hence, $\Pr[|X^*| \le 0.8h]$ is exponentially small in $h$.
\qed

We next need to bound the number of elements of $C^*(u,\F)$ with both endpoints in $X^*$.
\begin{lemma}\label{l:num-inc-x*}
With probability larger than $3/4$, $|C^*(u,\F) \cap (X^* \times X^*)| \le 1000 h^{2-\epsilon/2}$.
\end{lemma}
\Proof
We will actually upper bound the size of the potentially larger set $C^*(u,\F) \cap (X \times X)$.
Let $x,y$ be a pair of vertices in $C^*(u,\F)$.
If $x \in B(u,\E)$, then trivially $x \notin X$.
Similarly, if $y \in B(u,\E)$, then $y \notin X$.
Hence, we can assume that $x \in S'_i$ and $y \in S'_j$. Now if $i=j$ then at least one of them is not selected to $X$,
since from each part we selected at most one vertex to $X$. So we can assume $i \neq j$.
If $i \notin J$ then $x \notin X$ since we selected no vertex of $S'_i$ to $X$.
Likewise, if $j \notin J$ then $y \notin X$ since we selected no vertex of $S'_j$ to $X$.
So we can assume that $i \neq j$, $i,j \in J$. But recall that $J \subseteq K(u,\E)$, so $|S_i'| \ge n/(15h)$ and  $|S_j'| \ge n/(15h)$
(i.e. they are not very small). Since the vertex from $S'_i$ was chosen at random and the vertex of $S'_j$ was chosen at random,
the probability of the pair $x,y$ to be in $X \times X$ is at most $(15h/n)^2$.
By Lemma \ref{l:num-inconsistent}, $|C^*(u,\F)| \le n^2 h^{-\epsilon/2}$.
Hence the expected value of $|C^*(u,\F) \cap (X \times X)|$ is at most 
$$
n^2 h^{-\epsilon/2} \cdot \left(\frac{15h}{n}\right)^2 \le 225 h^{2-\epsilon/2}\;.
$$
By Markov's inequality the probability of being larger than four times the expectation is less than $1/4$.
\qed

We now want to show that the vertices of $X^*$ are not too concentrated among some parts of $P_1',\ldots,P_h'$,
namely that no part $P_i'$ contains too many vertices of $X^*$.
\begin{lemma}\label{l:non-concentration}
With probability larger than $3/4$, for all $i =1,\ldots,h$ we have $|P_i' \cap X^*| \le h/100$.
\end{lemma}
\Proof
We will actually prove the slightly stronger statement bounding $|P_i' \cap X|$.
Consider some part $P'_i$ and the disjoint parts $P'_{i,j} = P'_i \cap S'_j$ for $j \in J$
(notice that these disjoint parts might not cover all the vertices of $P'_i$, but we don't care about these
uncovered vertices since they cannot be selected to $X$).
Now, Since we selected one vertex from each $S'_j$, we have that $|P'_{i,j} \cap X| \in \{0,1\}$.
It is therefore an indicator random variable with probability of success
$$
\Pr[|P'_{i,j} \cap X| = 1] = \frac{|P'_{i,j}|}{|S'_j|} \le \frac{|P'_{i,j}|}{n/(15h)}
$$
where we have used the fact that $S'_j$ is not very small.
It now follows that the expectation of $|P_i' \cap X|$ is at most
$$
\sum_{j \in J} \frac{|P'_{i,j}|}{n/(15h)} \le \frac{15h}{n}|P_i'| \le \frac{15h}{n} \cdot 0.0004n \le 0.006h
$$
where we have used the fact that no single set $P_i$ (moreover $P_i' \subseteq P_i$) can contain more than $0.0004n$ vertices, by Lemma \ref{l:some-parts}.
But now observe that $|P_i' \cap X|$ is the sum of $|J|$ independent indicator random variables (recall, from each $S'_j$ for $j \in J$ we select a single
vertex to $X$ independently), so by a Chernoff bound (see appendix in \cite{AS-2004}), the probability that $|P_i' \cap X|$ is larger than its
expectation by a constant factor is exponentially small in $|J| \ge 0.9h$, so in particular, exponentially small in $h$.
So, for $h$ sufficiently large, with probability smaller than $1/(4h)$, $|P_i' \cap X| > h/100$.
Now, by the union bound, we have that with probability larger than $3/4$, for all $i=1,\ldots,h$ it holds that
$|P_i' \cap X| \le h/100$.
\qed

Since each of the four lemmas \ref{l:far-x}, \ref{l:x*}, \ref{l:num-inc-x*} \ref{l:non-concentration} 
states that the corresponding event happens with probability larger than $3/4$, we have that with positive probability, all four events hold.
Thus, we can pick a set $X=\{x_j~|~j \in J\}$ such that $X$ and its subset $X^* = X \cap W(u,\F)$  satisfy the following four properties:
\begin{enumerate}
\item[P1]
$G[X]$ is $4h^{2-\epsilon/2}$-close to $H[J]$. In particular, $G[X^*]$ is $4h^{2-\epsilon/2}$-close to $H[J^*]$ where $J^* \subset J$
satisfies $J^*= \{j~|~x_j \in X^*\}$.
\item[P2]
$|X^*|=|J^*| > 0.8h$.
\item[P3]
$|C^*(u,\F) \cap (X^* \times X^*)| \le 1000 h^{2-\epsilon/2}$.
\item[P4]
For all $i =1,\ldots,h$ we have $|P_i' \cap X^*| \le h/100$.
\end{enumerate}

We now consider the way the vertices of $X^*$ appear in $P_1',\ldots,P_h'$.
Let $S \subset [h] = \{i~|~X^* \cap P'_i \neq \emptyset\}$.
So $i \in S$ means that it contains a least one vertex of $X^*$ and $i \notin S$ means that it does not contain vertices of $X^*$.
There are two cases to consider:

\noindent
{\bf Case 1:} $|S| \le 0.7h$\;.
In  this case we obtain by P1 on the one hand, that $G[X^*]$ is $4h^{2-\epsilon/2}$-close to $H[J^*]$.
On the other hand, we obtain by P3 that there is a blowup of $H[S]$ that is blowup $1000 h^{2-\epsilon/2}$-close to $G[X^*]$.
Furthermore, $|S| \le 0.7h$ while $|X^*| \ge 0.8h$ by P2 and no blowup part contains more than $h/100$ vertices by P4.
So, this blowup is blowup $(4h^{2-\epsilon/2}+1000 h^{2-\epsilon/2})$-close to $H[J^*]$.
Since $4h^{2-\epsilon/2}+1000 h^{2-\epsilon/2} < 10^{-5}h^2$ this contradicts the fourth condition in the definition
of ${\cal P}_h$, thus our assumption that $H \in {\cal P}_h$.

\noindent
{\bf Case 2:} $|S| \ge 0.7h$\;.
Recall that we assume that $u \in P_i$ and $u \in S_j$ where $i \neq j$. So we fix the two indices $i$ and $j$.
Let $S^* \subseteq S \cup \{i\}$ have precisely $\lceil 0.7h \rceil$ elements where we force $i \in S^*$.
Now consider the following set of vertices $X^{**}$ of $G$. For each $t \in S^*$, if $t \neq i$ then add to $X^{**}$ some vertex of $X^*$
which belongs to $P_t'$ (this is always possible since $P'_t \cap X^* \neq \emptyset$ by the definition of $S$). If $t=i$, then add $u$ to $X^{**}$.
Observe that $|X^{**}|=|S^*|=\lceil 0.7h \rceil$ and that $u$ is always in $X^{**}$.
We will prove that $G[X^{**}]$ is very close to $H[S^*]$ and that it is also very close to
$H[Y^*]$ where $Y^* = \{ p ~|~ x_p \in X^{**}\} \cup \{j\}$. Observe also that $|Y^*|=\lceil 0.7h \rceil$.

Consider the natural bijection $\pi$ from $S^*$ to $X^{**}$ defined as follows. Let $t \in S^*$. If $t \neq i$ let $\pi(t)$ be the vertex of $X^{**}$ which belongs to
$P_t'$. For $t=i$ define $\pi(i)=u$.
Why is $\pi$ close to an isomorphism between $H[S^*]$ and $G[X^{**}]$? Every inconsistent pair has to be in $C^*(u,\F)$ and also has to be in $(X^* \times X^*)$, unless that pair involves $u$. But by the definition of $\F_u$, all the vertices in the parts $P_1',\ldots,P_h'$ are consistent with $u$.
Hence, by P3, $\pi$ is $1000 h^{2-\epsilon/2}$-close to an isomorphism between $H[S^*]$ and $G[X^{**}]$.

Consider the natural bijection $\sigma$ from $X^{**}$ to $Y^*$ defined as follows. Suppose $x_p \in X^{**}$ then $\sigma(x_p)=p$.
Also define $\sigma(u)=j$.
Why is $\pi$ close to an isomorphism between $G[X^{**}]$ and $H[Y^*]$? 
Every inconsistent pair that does not involve $u$, has to be in $G[X]$.
Also, $Y^* \setminus \{j\} \subset J^*$. But by the definition of $\E_u$, all the vertices in the parts $S_1',\ldots,S_h'$ are consistent with $u$.
Hence, by P1, $\sigma$ is $4h^{2-\epsilon/2}$-close to an isomorphism between $G[X^{**}]$ and $H[Y^*]$.

It follows that the composition $\pi\sigma$ is $4h^{2-\epsilon/2}+1000 h^{2-\epsilon/2}$ close to an isomorphism between  $H[S^*]$ and $H[Y^*]$.
Since $4h^{2-\epsilon/2}+1000 h^{2-\epsilon/2} < 10^{-5}h^2$, we must have by the fifth condition in the definition
of ${\cal P}_h$, that $\pi\sigma$ has less than $0.1h$ non-stationary points.

There now two options to consider. If $|S^* \; \triangle \; Y^*| \ge 0.2 h$, then trivially $\pi\sigma$ has at least $|S^* \triangle Y^*|/2 = 0.1h$
non-stationary points, a contradiction. If $|S^* \; \triangle \; Y^*| < 0.2 h$, then $|S^* \cap Y^*| \ge 0.6h$.
But since $\pi(i)=u$ and $\sigma(u)=j$ this means that in $H$, the distinct vertices $i$ and $j$ have more than $0.6h-2 > 0.55h$ vertices in their
agreement set, contradicting the assumption that $H \in {\cal P}_h$.
\qed

\section{Inducibility}

In this section we prove Theorem \ref{t:inducibility}.
Let $H \in {\cal P}_h$ from Theorem \ref{t:main}.
By Theorem \ref{t:main} it holds for $n \le 2^{\sqrt{h}}$ that $i_H(n)=g(n,h)$.
It will be profitable to use the largest $k$ such that $n=h^k$ and still $n \le 2^{\sqrt{h}}$.
We can therefore assume that $k \ge h^{1/3}+2$ as
$$
h^{\lceil h^{1/3}+2 \rceil} \le 2^{\sqrt{h}}\;.
$$
Now, by (\ref{e:ghk}) we have that
$$
i_H(h^k) = g(h^k,h) = \frac{h^{h(k-1)}(1- h^{k(1-h)})}{1-h^{1-h}}\;.
$$
Since $g(h^k,h)/\binom{h^k}{h}$ serves as a trivial upper bound for $i_H$ we obtain that
$$
i_H \le \frac{h^{h(k-1)}(1- h^{k(1-h)})}{(1-h^{1-h})\binom{h^k}{h}}\;.
$$
Clearly if we would have let $k \rightarrow \infty$ we would have obtained
$i_H = h!/(h^h-h)$. However we can only assume that $k$ is bounded from below $h^{1/3}+2$.
Nevertheless, we have that
$$
\frac{h^{h(k-1)}(1- h^{k(1-h)})}{(1-h^{1-h})\binom{h^k}{h}} \le 
\frac{h! h^{h(k-1)}(1- h^{k(1-h)})}{(1-h^{1-h})(h^k-h)^h} \le
\frac{h!}{h^h-h} \cdot \frac{1}{(1-h^{1-k})^h}\;.
$$
But notice that
$$
\frac{1}{(1-h^{1-k})^h} \le  e^{\frac{2}{h^{k-2}}} \le 1 + \frac{4}{h^{k-2}}
\le 1+\frac{4}{h^{h^{1/3}}}\;.
$$
It follows that
$$
i_H \le \frac{h!}{h^h-h} \cdot \left( 1+\frac{4}{h^{h^{1/3}}} \right)\;.
$$
\qed

\section{Concluding remarks and open problems}

We conjecture that the statement of Theorem \ref{t:super} and hence of Theorem \ref{t:main} can be extended to all $n$. Moreover:
\begin{conj}\label{c:1}
For all strongly asymmetric graphs $H$ it holds that $i_H(n)=g(n,h)$ for all $n \in \mathbb{N}$.
In particular,
$$
i_H= \frac{h!}{h^h-h}\;.
$$
\end{conj}
%Though seemingly easier, it would also be nice to prove that for all sufficiently large $h$, a random graph
%$H \sim G(h,1/2)$ has $i_H(n)=g(n,h)$ for all $n \in \mathbb{N}$ with probability $1-o_h(1)$.

It would be extremely interesting to determine all graphs for which $i_H = \frac{h!}{h^h-h}$.
As mentioned in the introduction, all cycles of length at least $5$ are conjectured to be in this family.
Likewise, determining the set of graphs for which $i_H(n)=g(n,h)$ for, say, all $n \le h^3$,
would also be gratifying as such graphs have $i_H= (1+o_h(1))\frac{h!}{h^h-h}$, as can be seen from the
proof of Theorem \ref{t:inducibility}. Given that this is true (for a much larger range of $n$) for all strongly asymmetric graphs,
one might stipulate that Theorem \ref{t:main} (or maybe even Conjecture \ref{c:1}) holds for all asymmetric graphs.
In fact, as the next proposition shows, this is very far from true.
\begin{prop}
For all $h$ sufficiently large, there are asymmetric graphs $H$ on $h$ vertices
for which $i_H(n) > g(n,h)$ for $n < h \log h$.
\end{prop}
\Proof
Recall that a random graph from $G(q,1/2)$ is almost surely asymmetric \cite{ER-1963}. Furthermore, almost surely a random graph from $G(q,1/2)$
has all its degrees at most $0.55q$ and the number of common neighbors of every two vertices is at least $0.2q$
(see the proofs of Lemma \ref{l:cond-1} and Lemma \ref{l:cond-2}) .
Hence, for $q$ sufficiently large, we can find (much more than) $2^q$ pairwise non-isomorphic graphs
$H_1,\ldots,H_{2^q}$ such that each $H_i$ is an asymmetric graph on $q$ vertices,
its maximum degree is at most $0.55q$, and the number of common neighbors in $H_i$ of every two vertices is at least $0.2q$.

Now let $H$ be the graph on $h=q2^q$ vertices obtained by taking a balanced blowup of the clique $K_{2^q}$
where each blowup part is of size $q$. Suppose the parts are $A_1,\ldots,A_{2^q}$.
Now make each $A_i$ induce a copy of $H_i$.

We claim that $H$ is asymmetric. Indeed, consider some automorphism $\pi$.
We first notice that $\pi(A_i) = A_j$ for some $j$ (possibly $i=j$).
Suppose not, then there are $u,v \in A_i$ such that $\pi(u) \in A_j$ and $\pi(v) \in A_k$ for $j \neq k$.
Since $\pi$ is an automorphism, the number of common neighbors of $u$ and $v$ must equal the number of common neighbors of
$\pi(u)$ and $\pi(v)$. But the number of common neighbors of $u$ and $v$ in $H$ is at least $h-q+z \ge h-q+0.2q$
where $z$ is the number of common neighbors of $u$ and $v$ in $H_i$.
On the other hand, the number of common neighbors of $\pi(u)$ and $\pi(v)$ in $H$ is $h-2q+x+y \le h-2q+0.55q+0.55q=h-q+0.1q$
where $x$ is the degree of $\pi(v)$ in $H_k$ and $y$ is the degree of $\pi(u)$ in $H_j$.
This contradicts that $\pi$ is an automorphism. Hence, we must have $\pi(A_i)=A_j$ for some $j$.

But now we claim that we must have $\pi(A_i)=A_i$. Indeed, if $\pi(A_i)=A_j$, then since $\pi$ is an automorphism, this means that
$H_i$ and $H_j$ are isomorphic. Since $H_i$ is not isomorphic to $H_j$ whenever $i \neq j$, this implies that $\pi(A_i)=A_i$
for all $i=1,\ldots, 2^q$. But if $\pi(A_i)=A_i$ this means that $\pi$ restricted to $A_i$ is an automorphism of $H_i$.
Since $H_i$ is asymmetric, $\pi$ must be trivial on $A_i$. As this holds for all $i=1,\ldots, 2^q$, we must have that $\pi$ is trivial, proving
that $H$ is asymmetric.

Designate a set $X \subset V(H)$ order $2^q$ obtained by taking one vertex from each $A_i$ for $i=1,\ldots,2^q$.
Notice in fact that $H[X]$ is a clique of order $2^q$.
Now consider the graph $G$ on $n=hq$ vertices which is the balanced blowup of $H$.
So the parts of this blowup are $\{P_x~|~x \in V(H)\}$ and $|P_x|=q$.
Now, replace each independent set $P_x$ with a copy of $H_i$ where $x \in A_i$ (recall that the $A_i$ are a partition of $H$ and each induce
a copy of $H_i$ in $H$). Now consider the subgraph of $G$ induced by the parts $\{P_x~|~x \in X\}$. Notice that this subgraph
is isomorphic to $H$. On the other hand, it is not an induced  copy of $H$ obtained by taking one vertex from each part of $G$,
so it is an additional induced copy to the $f(n,h)=g(n,h)=q^h$ induced copies that can be obtained as such.
Hence, $i_H(G) > g(n,h)$ while $n=hq < h \log h$.
\qed

Theorem \ref{t:super}, and therefore Theorem \ref{t:main} have variants which apply to other graph densities and some other combinatorial structures.
Indeed, suppose $0 < p < 1$ is given. One can adjust the definition of ${\cal P}_h$ to obtain a graph property ${\cal P}_h(p)$ which contains
almost all graphs of $G(h,p)$. Indeed, Condition 1 in the definition could be changed to requiring that the degree be between
$p(1-\delta)n$ and $p(1+\delta)n$ for an absolute constant $\delta$, that the agreement sets in Condition 2 would be
at most $(1+\delta)(p^2+(1-p)^2)$ and $(1+\delta)(p^3+(1-p)^3)$ respectively, that the distinguishing set in Condition 3 be of size at
most $-3\log (p^2+(1-p)^2) \log h$, and that the constant $10^{-5}h^2$ in Conditions 4,5 be replaced with $c\min \{p,1-p\}h^2$ for some
absolute constant $c$. The proof of Theorem \ref{t:super} and its lemmas stay essentially the same, after adjusting constants
everywhere. Similarly, a variant of the theorems applies to the family of tournaments.
We can naturally define a tournament to be strongly asymmetric
by modifying the definition of ${\cal P}_h$ to suit tournaments (namely, degree requirements are replaced with in-degree and out-degree requirements,
changing an edge means flipping its direction, etc.). Likewise, instead of considering $G(h,1/2)$ we consider the probability space of all
tournaments on $h$ vertices.

\section{Acknowledgment}
Independently of our work, Jacob Fox, Hao Huang, and Choongbum Lee (private communication) have sent the author a manuscript proving that almost all graphs $H$ 
have $i_H(n)=g(n,h)$ for all $n \in \mathbb{N}$. Once their paper is publicly available, a reference will be given.

\bibliographystyle{plain}

\bibliography{references}

\end{document}